\documentclass[runningheads,a4paper]{llncs}

\usepackage{amssymb}
\usepackage{graphicx}
\usepackage{amsmath,amsfonts,bm,mathtools}
\usepackage{algorithmicx}
\usepackage{algpseudocode}
\usepackage{booktabs}
\usepackage{threeparttable}
\usepackage{url}
\usepackage{tabularx}
\urldef{\mailsa}\path|{alfred.hofmann, ursula.barth, ingrid.haas, frank.holzwarth,|
\urldef{\mailsb}\path|anna.kramer, leonie.kunz, christine.reiss, nicole.sator,|
\urldef{\mailsc}\path|erika.siebert-cole, peter.strasser, lncs}@springer.com|
\newcommand{\keywords}[1]{\par\addvspace\baselineskip
\noindent\keywordname\enspace\ignorespaces#1}

\begin{document}

\mainmatter  

\title{A Multiobjective State Transition Algorithm for Single Machine Scheduling}

\titlerunning{A Multiobjective State Transition Algorithm for Single Machine Scheduling}

%
%
\author{Xiaojun Zhou$^{\dag}$%
\thanks{Corresponding author for this paper.
This research work is conducted between Deakin University and Ballarat University under the Collaboration Research Network (CRN) initiative. 
The problem studied in this paper is related to the Australian Research Council (ARC) linkage project number LP0991175.
}%
\and Samer Hanoun$^{\ddag}$
\and David Yang Gao$^{\dag}$
\and Saeid Nahavandi$^{\ddag}$}
%
\authorrunning{A Multiobjective State Transition Algorithm for Single Machine Scheduling}

\institute{$^{\dag}$School of Science, Information Technology and Engineering, University of Ballarat, Victoria 3353, Australia\\
$^{\ddag}$Centre for Intelligent Systems Research, Deakin University, Geelong, Australia}

%
%

\maketitle

\begin{abstract}
In this paper, a discrete state transition algorithm is introduced to solve a multiobjective single machine job shop scheduling problem. In the proposed approach, a non-dominated sort technique is used to select the best from a candidate state set, and a Pareto archived strategy is adopted to keep all the non-dominated solutions. Compared with the enumeration and other heuristics, experimental results have demonstrated the effectiveness of the multiobjective state transition algorithm.
\end{abstract}

\section{Introduction}
The multiobjective optimization is encountered in many real world applications \cite{Marler}. For a specific policy, the decision maker may find it advantageous for one goal but disadvantageous for others. A traditional way to deal with this issue is to impose a priori preference reflecting the relative importance of different objectives; however, the final solution just indicates a decision maker's satisfaction, and it might be dissatisfactory for other decision makers.

To ameliorate the problem, the concept of \textit{Pareto optimality} and other relevant concepts are introduced. These are defined as follows:

1) \textit{Pareto dominance}: A feasible solution $\bm x = (x_1, \cdots, x_n)$ is said to Pareto dominate another feasible solution $\bm y = (y_1, \cdots, y_n)$, denoted as $\bm x \prec \bm y$, if
\begin{eqnarray}
f_i(\bm x) \leq f_i(\bm y), \forall i \in \{1,\cdots,k\}, \;\mathrm{and}\;\; \exists j \in \{1,\cdots,k\}, f_j(\bm x) < f_j(\bm y),
\end{eqnarray}
where, $f_i(\bm x)$ is the $i$th objective function, $k$ is the number of objectives.

2) \textit{Pareto optimality}: A feasible solution $\bm x^{*}$ is said to be Pareto optimal if and only if
\begin{eqnarray}
\neg \exists \bm x \in S, \bm x \prec \bm x^{*},
\end{eqnarray}
where, $S$ is the feasible space.

3) \textit{Pareto optimal set}: The Pareto optimal set, denoted as $P^{*}$, is defined by
\begin{eqnarray}
P^{*} = \{\bm x^{*} \in S | \neg \exists \bm x \in S, \bm x \prec \bm x^{*}\}.
\end{eqnarray}

4) \textit{Pareto front}: The Pareto front, denoted as $P f^{*}$, is defined by
\begin{eqnarray}
P f^{*} = \{(f_1(\bm x^{*}), \cdots, f_k(\bm x^{*}))| \bm x^{*} \in P^{*}\}.
\end{eqnarray}

The introduction of \textit{Pareto optimality} allows us to find a set of Pareto optimal solutions simultaneously, independent of the decision maker's priori preference.

In the past few decades, evolutionary-based and nature-inspired multiobjective optimization techniques have drawn considerable attention for scheduling problems\cite{Coello,Lei,Al-Anzi,Xia,Hanoun2011,Hanoun2012}. In this paper, we introduce a recently new heuristics called state transition algorithm \cite{xzhou2011a,xzhou2011b,yang,xzhou2012} as the basic search engine for the multiobjective optimization. A non-dominated sort approach is used to select the best from a candidate state set, and the best state is stored using
a Pareto archive strategy. Experimental results have testified the effectiveness of the proposed algorithm.

\section{Problem Description}
In the field of joinery manufacturing, jobs with similar materials can be scheduled together to minimize the amount of materials used; therefore, reducing the cost.

For example, based on the cost savings matrix shown in Table \ref{costex1}, pairing Job1 and Job2 will provide saving in the cost equivalent to 4 units.

\begin{table}[!htbp]
\centering
\caption{The cost savings matrix for 5 jobs having the same material}
\label{costex1}
\begin{tabular}{p{1.5cm}p{1.5cm}p{1.5cm}p{1.5cm}p{1.5cm}p{1.5cm}}
\hline
& Job1 & Job2 & Job3 & Job4 & Job5 \\
\hline
Job1 & 0	& 4	    & 2.64	& 4.08  & 3.9 \\
Job2 & 4	& 0	    & 3.64  & 4.72	& 4.23\\
Job3 & 2.64 & 3.64	& 0	    & 2.65	& 2.87\\
Job4 & 4.08	& 4.72	& 2.65	& 0	    & 3.84\\
Job5 & 3.9	& 4.23	& 2.87	& 3.84	& 0\\
\hline
\end{tabular}
\end{table}

Additionally, based on the jobs' processing times and due dates as shown in Table \ref{dueex1}, and for any given sequence and pairing of jobs, not only the total cost saving C is affected but also the total tardiness time T, which is calculated as:
\begin{eqnarray}
T = \sum_{j=1}^{n}\max\{0,c_j - d_j\}
\end{eqnarray}
where, $c_j$ and $d_j$ are the completion time and the due time of job $j$, respectively.

\begin{table}[!htbp]
\centering
\begin{threeparttable}
\caption{Due dates and processing times for a set of 5 jobs}
\label{dueex1}
\begin{tabular}{p{2.5cm}p{2.5cm}p{4cm}}
    \hline
    Job & Due Date(days)& Processing Time(hours)\\
    \hline
    Job1 & 8 & 17:40\\
    Job2 & 2 & 24:00\\
    Job3 & 11 & 19:20\\
    Job4 & 3 & 25:00\\
    Job5 & 3 & 14:40\\
    \hline
\end{tabular}
\begin{tablenotes}
 \item [a]  Number of operational hours = 8 hours per day
\end{tablenotes}
\end{threeparttable}
\end{table}

The goal of this paper is to determine the optimal sequence with pairing, in order to maximize the total cost savings and minimize the total tardiness time.

It is obvious that finding the permutation of the sequence $\{1,2,\cdots,n\}$ with pairing becomes a solution to the multiobjective single machine scheduling problem; however, not without the necessity to discuss the number of pairs for any fixed sequence of jobs.

Given a sequence $s = (1,2,\cdots,n)$, for $n=3$, we have 2 possible pairing options (1-2)-3 and 1-(2-3); for $n=4$, we have 2 possible pairing options (1-2)-(3-4) and 1-(2-3)-4, as pairing options (1-2)-3-4 and 1-2-(3-4) are discarded; for $n=5$, we have 3 possible options (1-2)-(3-4)-5, (1-2)-3-(4-5) and 1-(2-3)-(4-5), as options  (1-2)-3-4-5, 1-2-(3-4)-5, 1-2-3-(4-5) and 1-(2-3)-4-5 are discarded.

If $P1(n)$ denotes the number of pairs with the first two jobs pairing, and $P2(n)$ denotes the complement of $P1(n)$, then we have the following theorem:
\\
\indent \textbf{Theorem 1}
\begin{eqnarray}
P1(n+1) = P(n-1), P2(n+1) = P1(n), n \geq 3
\end{eqnarray}
where, $P(n) = P1(n) + P2(n)$ is the total number of pairs.
For example, $P1(2) = 1, P2(2) = 0$, $P1(3) =1, P2(3) = 1$, $P1(4) = 1, P2(4) =1$, $P1(5) = 2, P2(5) = 1$, we have
$P1(4) = P(2), P2(4) = P1(3)$, $P1(5) = P(3), P2(5) = P1(4)$.

\begin{figure}[!htbp]
  \centering
  \includegraphics[width=8cm,height=6cm]{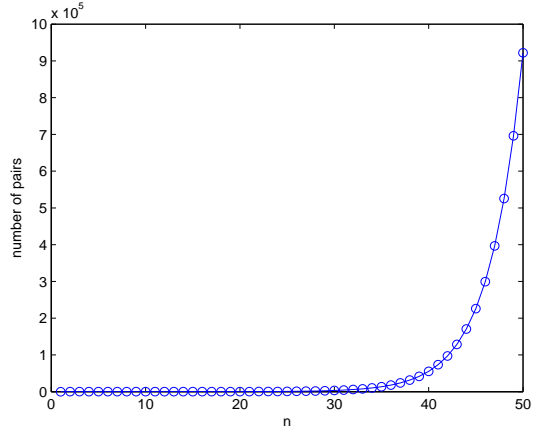}
  \caption{Growth trend relative to the sequence size}\label{trend}
\end{figure}

\bigskip
Fig. \ref{trend}, shows the growth trend of the number of pairs with the sequence size; however, only small size job scheduling problems are considered in this study.
Considering that $P(10) = 12 \ll 10! = 3628800$, a complete enumeration approach is used for pairing and only the permutation of a sequence is focused.
\section{Discrete State Transition Algorithm}
In the case a solution to a specific optimization problem is described as a state, then the transformation to update the solution becomes a state transition. Without loss of generality, the unified form of discrete state transition algorithm can be described as:
\begin{eqnarray}
\left \{ \begin{array}{ll}
\bm x_{k+1}= A_{k}(\bm x_{k}) \bigoplus B_{k}(\bm u_{k})\\
y_{k+1}= f(\bm x_{k+1})
\end{array} \right.,
\end{eqnarray}
where, $\bm x_{k} \in \mathcal{Z}^{n}$ stands for a current state, corresponding to a solution of a specific optimization problem; $\bm u_{k}$ is a function of $\bm x_{k}$ and historical states; $A_{k}(\cdot)$, $B_{k}(\cdot)$ are transformation operators, which are usually state transition matrixes; $\bigoplus$ is an operation, which is admissible to operate on two states; and $f$ is the cost function or evaluation function.

\begin{figure}[!htbp]
  \centering
  \includegraphics[width=8cm]{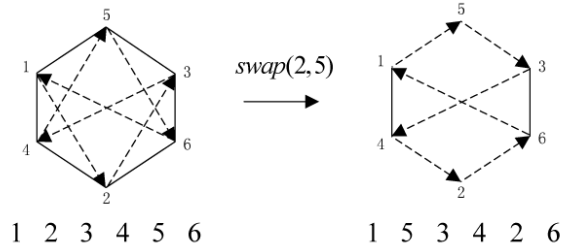}
  \caption{Illustration of the swap transformation}\label{swaptsp}
\end{figure}

\begin{figure}[!htbp]
  \centering
  \includegraphics[width=8cm]{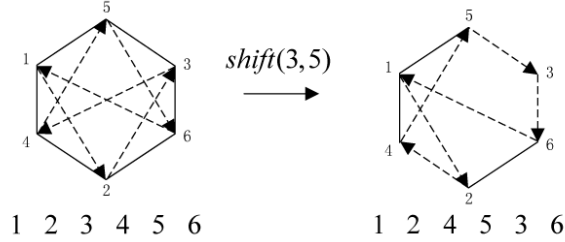}
  \caption{Illustration of the shift transformation}\label{shifttsp}
\end{figure}

\begin{figure}[!htbp]
  \centering
  \includegraphics[width=8cm]{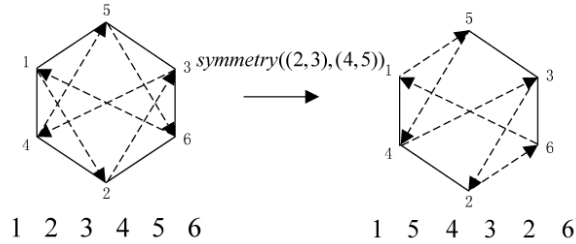}
  \caption{Illustration of symmetry transformation}\label{symtsp}
\end{figure}

The following three transformation operators are defined to permute current solution \cite{yang}:\\
(1) Swap Transformation
\begin{eqnarray}
\bm x_{k+1}= A^{swap}_{k}(m_a) \bm x_{k},
\end{eqnarray}
where, $A^{swap}_{k} \in \mathbb{R}^{n \times n}$ is the swap transformation matrix, $m_a$ is the swap factor, a constant integer used to control the maximum number of positions to be exchanged, while the positions are random. Fig. \ref{swaptsp} shows an example of the swap transformation with $m_a = 2$.
\\
(2) Shift Transformation
\begin{eqnarray}
\bm x_{k+1}= A^{shift}_{k}(m_b) \bm x_{k},
\end{eqnarray}
where, $A^{shift}_{k} \in \mathbb{R}^{n \times n}$ is the shift transformation matrix, $m_b$ is the shift factor, a constant integer used to control the maximum length of consecutive positions to be shifted. Note that both the selected position to be shifted after and positions to be shifted are chosen randomly. Fig. \ref{shifttsp} shows an example of the shift transformation with $m_b = 1$.
\\
(3) Symmetry Transformation
\begin{eqnarray}
\bm x_{k+1}= A^{sym}_{k}(m_c) \bm x_{k},
\end{eqnarray}
where, $A^{sym}_{k} \in \mathbb{R}^{n \times n}$ is the symmetry transformation matrix, $m_c$ is the symmetry factor, a constant integer used to control the maximum length of subsequent positions as center. Note that both the component before the subsequent positions and consecutive positions to be symmetrized are created randomly. Fig. \ref{symtsp} shows an example of the symmetry transformation with $m_c = 0$.

\section{Pareto Archived Strategy based on DSTA}
In state transition algorithm, the times of transformation are called search enforcement (\textit{SE}); as a result,
after each transformation operator, a candidate state set $S$ is generated.

\subsection{Non-dominated sort}
We use a sorting approach similar to the fast-non-dominated-sort proposed in \cite{Deb}, described as follows:
\begin{algorithmic}[1]
\For{each $s \in S$}
    \State $n_s \gets 0$
    \For{each $t \in S$}
        \If{$t \prec s$}
             \State $n_s \gets n_s + 1$
        \EndIf
    \EndFor
\EndFor
\end{algorithmic}
where, $n_s$ is the domination count, representing the number of solutions dominating solution $s$.
After the non-dominated sort, the state with the least count will be stored as incumbent \textit{best} for the next transformation operator.

\subsection{Pareto Archived Strategy}
We adopt a simple Pareto archived strategy to select current \textit{best}, as follows:
\begin{algorithmic}[1]
\For{each $A_i \in \mathcal{A}$}
    \If{$best \prec A_i$}
         \State $\mathcal{A} \gets \mathcal{A} - A_i$
    \ElsIf{$A_i \prec best$}
         \State $\mathcal{A} \gets \mathcal{A} $
    \Else
         \State $\mathcal{A} \gets \mathcal{A} \bigcup best$
    \EndIf
\EndFor
\end{algorithmic}
where, $\mathcal{A}$ is the archive keeping all non-donominated solutions.
\subsection{Pseudocodes of the proposed algorithm}
The core procedure of the proposed algorithm can be outlined in pseudocodes:
\begin{algorithmic}[1]
\Repeat
    \State $State \gets operator(best,SE,n)$
    \State $best \gets update\_best(best,SE,n,data)$
    \State $Paretoset \gets update\_archive(Paretoset,best)$
\Until{the maximum number of iterations is met}
\end{algorithmic}
where, $State$ is the state set;  $operator$ stands for the three transformation operators, which are carried out sequentially;  $update\_best$ is corresponding to the
non-dominated sort, and $update\_archive$ corresponds to the Pareto archived strategy. The $data$ is the known information (cost saving matrix, due dates and processing times) about a specific scheduling problem.
\section{Experimental Results}
In order to test the performance of the proposed multiobjective state transition algorithm, two typical examples are used for comparison.
In the following experiments, $SE = 20, m_a = 2, m_b = 1, m_c = 0$ are adopted for parameter settings. The maximum number of iterations for
are 100 and 1000 respectively for the two examples.
The known data for Example 1, 2 are given in Table \ref{costex1} and Table \ref{dueex1}, Table \ref{costex2} and Table \ref{dueex2}, respectively, and the corresponding results can be found in Table \ref{resultex1} and Table \ref{resultex2}. It is worth to note that, the pairing methodology used with complete enumeration and Cuckoo Search (CS) is based on a greedy approach by first selecting the pair that produces the highest cost savings, and then repeating the same procedure for the remaining set of pairs in the sequence \cite{Hanoun2012}. We can find that for Example 1, STA obtained a solution which can dominate the optimal solutions by enumeration and CS. From both examples, it is easy to find that some additional optimal solutions are achieved by STA.
\begin{table}[!htbp]
\centering
\caption{The cost savings matrix for 10 jobs having the same material}
\label{costex2}
\begin{tabular}{p{1cm}p{1cm}p{1cm}p{1cm}p{1cm}p{1cm}p{1cm}p{1cm}p{1cm}p{1cm}p{1cm}}
  \hline
  & Job1 & Job2 & Job3 & Job4 & Job5 & Job6 & Job7 & Job8 & Job9 & Job10\\
  \hline
Job1&0	& 2.73	& 2.1&	2.16&	2.66&	3.6	&2.46&	2.7	&2.46&	2.8\\
Job2&2.73&	0	&2&	1.6	&4.3&	3.69&	2.3&	3.5	&2.76&	3.6\\
Job3&  2.1&	2	&0	&1.4&	3.51	&3.33&	2.52	&3.68	&2.52&	2.46\\
Job4&  2.16	&1.6&	1.4	&0	&2.17	&2.32&	2.72	&3.04&	2.04&	2.97\\
Job5&  2.66&	4.3	&3.51	&2.17&	0	&3.6	&4.05&	4.41	&2.7&	2.64\\
Job6&  3.6	&3.69&	3.33	&2.32&	3.6&	0&	2.58	&4.7	&3.44&	2.94\\
Job7&  2.46&	2.3&	2.52&	2.72	&4.05	&2.58	&0	&2.6	&2.88&	2.82\\
Job8&  2.7	&3.5&	3.68	&3.04	&4.41	&4.7	&2.6	&0&	3.64	&3.57\\
Job9&  2.46	&2.76	&2.52&	2.04&	2.7	&3.44	&2.88&	3.64&	0&	3.76\\
Job10& 2.8	&3.6	&2.46	&2.97	& 2.64&	2.94&	2.82	&3.57&	3.76&	0\\
\hline
\end{tabular}
\end{table}

\begin{table}[!htbp]
\centering
\caption{Due dates and processing times for a set of 10 jobs}
\label{dueex2}
\begin{tabular}{p{2.5cm}p{2.5cm}p{4cm}}
    \hline
    Job & Due Date(days)& Processing Time(hours)\\
    \hline
    Job1 & 11 & 14:00\\
    Job2 & 2 & 18:00\\
    Job3 & 13 & 15:00\\
    Job4 & 14  & 8:20\\
    Job5 & 11 & 17:20\\
    Job6 & 9 & 16:00\\
    Job7 & 4 & 19:40\\
    Job8 & 6 & 23:20\\
    Job9 & 10 & 20:00\\
    Job10& 10 & 19:20\\
    \hline
\end{tabular}
\end{table}

\begin{table}[!htbp]
\centering
\caption{Comparison results for the set of jobs presented in Table \ref{costex1} and Table \ref{dueex1}}
\label{resultex1}
\begin{tabular}{p{2cm}p{4cm}p{1cm}p{1cm}}
  \hline
  Approach & Optimal solutions & T & C \\
  \hline
   Complete    & (2-5)-(1-4)-3 & 13 & 8.31\\
   Enumeration & (5-2)-(1-4)-3 & 13 & 8.31\\
               & (2-5)-(4-1)-3 & 13 & 8.31\\
               & (2-4)-(5-1)-3 & 15 & 8.62\\
  \hline
  CS\cite{Hanoun2012}  & (2-5)-(1-4)-3 & 13 & 8.31\\
                       & (5-2)-(1-4)-3 & 13 & 8.31\\
                       & (2-5)-(4-1)-3 & 13 & 8.31\\
                       & (2-4)-(5-1)-3 & 15 & 8.62\\
  \hline
  STA         & (5-2)-(1-4)-3 & 13 & 8.31\\
              & (2-5)-(1-4)-3 & 13 & 8.31\\
              & (2-5)-(4-1)-3 & 13 & 8.31\\
              & (5-2)-(4-1)-3 & 13 & 8.31\\
              & \textbf{(2-4)-(5-1)-3} & \textbf{15} & \textbf{8.62}\\
  \hline
\end{tabular}
\end{table}

\begin{table}[!htbp]
\centering
\caption{Comparison results for the set of jobs presented in Table \ref{costex2} and Table \ref{dueex2}}
\label{resultex2}
\begin{tabular}{p{2cm}p{6cm}p{1cm}p{1cm}}
  \hline
  Approach & Optimal solutions & T & C \\
  \hline
Complete    & (5-7)-(2-6)-(1-3)-(4-10)-(8-9) & 39 &16.45\\
Enumeration & (5-7)-(2-6)-(1-3)-(4-8)-(10-9)& 40 & 16.64\\
            & (5-7)-(2-6)-(1-3)-(4-8)-(9-10)& 40 & 16.64\\
            & (5-7)-(2-6)-(1-4)-(3-8)-(10-9)& 41 & 17.34\\
            & (5-7)-(2-6)-(1-4)-(3-8)-(9-10)& 41 & 17.34\\
            & (5-2)-(7-4)-(6-1)-(3-8)-(10-9)& 43 & 18.06\\
            & (5-2)-(7-4)-(6-1)-(3-8)-(9-10)& 43 & 18.06\\
            & (2-5)-(7-4)-(6-1)-(3-8)-(10-9)& 43 & 18.06\\
            & (2-5)-(7-4)-(6-1)-(3-8)-(9-10)& 43 & 18.06\\
  \hline
  CS\cite{Hanoun2012} & 2-(7-5)-(6-1)-3-(4-10)-(8-9) & 39 & 14.26\\
                      & (5-7)-(2-6)-(1-3)-(4-8)-(9-10)& 40 & 16.64\\
                      & (5-7)-(2-6)-(1-4)-(3-8)-(10-9)& 41 & 17.34\\
                      & (2-5)-(7-4)-(6-1)-(3-8)-(10-9)& 43 & 18.06\\
                      & (2-5)-(7-4)-(6-1)-(3-8)-(9-10)& 43 & 18.06\\
                      & (5-2)-(7-4)-(6-1)-(3-8)-(10-9)& 43 & 18.06\\
  \hline
  STA         & (5-7)-(2-6)-(1-3)-(4-10)-(8-9) & 39 & 16.45\\
              & \textbf{(5-7)-(2-6)-(1-3)-(4-10)-(9-8)} & \textbf{39} & \textbf{16.45}\\
              & (5-7)-(2-6)-(1-3)-(4-8)-(10-9) & 40 & 16.64\\
              & (5-7)-(2-6)-(1-3)-(4-8)-(9-10) & 40 & 16.64\\
              & (5-7)-(2-6)-(1-4)-(3-8)-(10-9) & 41 & 17.34\\
              & (5-7)-(2-6)-(1-4)-(3-8)-(9-10) & 41 & 17.34\\
              & (5-2)-(7-4)-(6-1)-(3-8)-(10-9) & 43 & 18.06\\
              & (5-2)-(7-4)-(6-1)-(3-8)-(9-10) & 43 & 18.06\\
              & (2-5)-(7-4)-(6-1)-(3-8)-(10-9) & 43 & 18.06\\
              & (2-5)-(7-4)-(6-1)-(3-8)-(9-10) & 43 & 18.06\\
  \hline
\end{tabular}
\end{table}
\section{Conclusion}
A multiobjective state transition algorithm is presented for a single machine job shop scheduling problem. In this paper, a complete enumeration approach is used for
pairing the jobs in a fixed sequence. Compared with a greedy-based approach used with both the complete enumeration method and the CS, experimental results show the effectiveness of the proposed algorithm in obtaining the true set of all Pareto optimal solutions.


\end{document}